\documentclass[12pt]{amsart}
  \usepackage{amssymb}
 \newcommand{\bee}[1]{\begin{equation}\label{#1}}
 \newcommand{\ene}{\end{equation}}

\begin{document}

MSC2010 17D99

\medskip
\begin{center}
{\Large \bf A new properties of varieties of Leibnitz algebras}

\medskip
{\large A. V. Shvetsova, T. V. Skoraya}
\end{center}

The paper is devoted to the study of the new properties of varieties
of Leibnitz algebras. The characteristic of base field $\Phi$
assumed to be zero. All undefined concepts can be found in
\cite{GA-ZMV}. The article presented two new results. The first
result belongs to the second author and contains a proof of the
sufficient conditions for the finiteness of the colength varieties
of Leibnitz algebras. The second belongs to the first author. In it
are found a basis of identities and a basis of the space of
multilinear elements of variety $\widetilde{{\bf V}}_{3}$ of
Leibnitz algebras.

\section*{1. Introduction}

A linear algebra with bilinear multiplication, which is satisfies to
the Leibnitz identity $ (xy)z \equiv (xz)y+x(yz),$ is called a
Leibnitz algebra. Perhaps for the first time this concept was
discussed in the article \cite{BAM} as a generalization of Lie
algebras. The Leibnitz identity allows any element expressed as a
linear combination of elements in which the brackets are arranged
from left to right. Therefore further agree omit brackets in
left-normed products, i.e. $(((ab)c)\dots d) = abc\dots d.$ A
variety ${\bf V}$ of linear algebras over a field $\Phi$ is a set of
algebras over this field that satisfy a fixed set of identities.
Note, that the system of identities can be given implicitly. In this
case the variety $\textbf{V}$  is usually defined generating algebra
given constructively.

Let $F(X,{\bf V})$ be a relatively free algebra of variety ${\bf V}$
with countable set of free generators $X = \{x_1,x_2,\dots\}$.
Consider the space of the multilinear elements of algebra $F(X,{\bf
V})$. This space we will denote $P_n({\bf V})$ and call multilinear
part of variety ${\bf V}$. On this space naturally introduce the
action of permutations that we can consider it as a $\Phi
S_n$-module, where $S_n$ is a symmetric group. Since the field
$\Phi$ has zero characteristic, then the space $P_n({\bf V})$ is the
direct sum of irreducible submodules. Denote by $\chi_\lambda$ the
character of the irreducible representations of the symmetric group,
which corresponds to the partition $\lambda$ of the number $n$. Then
the character of module $P_n({\bf V})$ is expressed by the formula
\begin {equation}
\label {haracter1} \chi(P_n({\bf V}))=\sum_{\lambda\vdash n}
m_\lambda \chi_\lambda,
\end {equation}
where $m_\lambda$ are the multiplicities of irreducible submodules
in this sum.

An important numerical characteristic of variety ${\bf V}$ of linear
algebras is the colength $l_n({\bf V})$, which is defined as the
number of terms in the decomposition of character in the sum of
irreducible characters:
\begin {equation}
\label {colength2} l_n({\bf V}) = \sum_{\lambda\vdash n} m_\lambda.
\end {equation} We say that the colength ${\bf V}$ of variety is finite if
there exists a constant $C$ independent of $n$ such that for any $n$
 is true the inequality $l_n({\bf V}) \leq C.$

The right multiplication operator, for example, on element $z$, we
denote by $Z$, assuming that $xz=xZ$. This designation allows the
element $x {\underbrace{y...y}_n}$ to write in the form $xY^n$.
Recall that the standard polynomial of degree $n$ has the form:
$St_{n}(x_{1},x_{2},\dots ,x_{n})=\sum_{q\in
S_{n}}(-1)^{q}x_{q(1)}x_{q(2)}\dots x_{q(n)},$ where the summation
is carried out by elements of the symmetric group, and $(-1)^{q}$ is
equal to $+1$ or $-1$ depending on the parity of permutation $q$.
Agree variables in standard polynomial denote with special
characters above (below, wave and so on). For example the standard
polynomial of degree $n$ in the variables $x_{1},x_{2},\dots,$
$x_{n}$ we will write as follows:
$St_{n}=\overline{x}_{1}\overline{x}_{2}\dots\overline{x}_{n}$. It
is clear that the standard polynomial is skew symmetric. Variables
in different skew symmetric sets will be denoted by different
symbols, for example: $$\sum_{q\in S_{n}, p\in
S_{m}}(-1)^{q}(-1)^{p}x_{q(1)}x_{q(2)}\dots
x_{q(n)}y_{p(1)}y_{p(2)}\dots y_{p(m)}=$$
$$=\overline{x}_{1}\overline{x}_{2}\dots
\overline{x}_{n}\widetilde{y}_{1}\widetilde{y}_{2}\dots
\widetilde{y}_{m}.$$

\section*{2. Sufficient condition for the finiteness colength}
\section*{of varieties of Leibnitz algebras}

Previously, in the article \cite{SAV} are identified the necessary
conditions for the finiteness colength of varieties of Leibnitz
algebras. Further, we consider it sufficient conditions.

Following article \cite{RSM} denote the variety of all Leibnitz
algebras (Lie algebras), that satisfy the identity
$(x_1x_2)(x_3x_4)\dots (x_{2s+1}x_{2s+2}) \equiv 0$, by
$\widetilde{{\bf N}_s{\bf A}}$ (relatively ${\bf N}_s{\bf A}$). Let
in addition ${\bf V}_{1}={\bf N}_2{\bf A}$ is a variety of all Lie
algebras, commutator of which is nilpotent of class not more then
two, and $\widetilde{{\bf V}}_{1}$ is a variety of Leibnitz algebras
defined by the identity $x_{1}(x_{2}x_{3})(x_{4}x_{5})\equiv 0$.

\textbf{Theorem 1.} \textit{Let ${\bf V}$ be a subvariety of variety
 $\widetilde{{\bf N}_s{\bf A}},$ which for any natural numbers $k, m,\, k\le m,$ and $\alpha_1,
\dots , \alpha_k \in K$ satisfies the identity
\begin {equation}
\label {th3} xY^kzY^{m-k} \equiv \sum_{i=1}^k \alpha_i
xY^{k-i}zY^{m-k+i}.
\end {equation}
Then the variety ${\bf V}$ has the final colength.}

\textbf{Proof.} Because identity (\ref{th3}) is not satisfied in the
varieties ${\bf V}_{1}$ and $\widetilde{{\bf V}}_{1}$, from
conditions of the theorem follows that ${\bf V_{1}}, \widetilde
{{\bf V}}_{1}\not\subset {\bf V} \subset {\widetilde{\bf N_sA}}$.
Then by theorem 1 of article \cite{MSP-COI}, there exists a constant
$C$ independent of $n$ for such in the sum (\ref{haracter1}) is true
the condition $(n-\lambda_{1})<C$. In this case in the sum
(\ref{colength2}) the number of the non-zero terms is bounded by a
constant independent of $n$. Thus, to prove the result, it suffices
to establish that all multiplicities $m_{\lambda}$ are bounded by a
constant, which also is independent of $n$.

The article \cite{MSP-ZMV} is proved that the multiplicity
$m_{\lambda}({\bf V})$ is equal to the number of linearly
independent polyhomogeneous elements of special form. We will show
that the dimension of the space of polyhomogeneous elements is
bounded by a constant independent of $n$, which will complete the
proof of the theorem.

Consider $\lambda$, for which $m_{\lambda}\neq 0$. For such
partition is true the condition $(n-\lambda_{1})<C$ and to it
correspond monomials of the form
$$g_s=Y^{\alpha_{1}}x_{i_{1}}Y^{\alpha_{2}}x_{i_{2}}Y^{\alpha_{3}}\dots
Y^{\alpha_{s}}x_{i_{s}}Y^{\alpha_{s+1}},$$ where $s<C$. Denote by
$Q_{\lambda_{1}}$ the space generated by elements $g_{s}$. We prove
that the number of linearly independent monomials $g_{s}$ bounded by
a constant. The proof is by induction on the number $s$ of
generators $x_{i}$ and lexicographic order on lines of the form
$(\alpha_{1}, \alpha_{2},\dots, \alpha_{s+1})$.

Consider the case $s=1$. Then generating monomials of the space
$Q_{\lambda_{1}}$ have the form:
$Y^{\alpha_{1}}x_{1}Y^{\alpha_{2}}.$ If for these elements are true
the conditions $\alpha_{1}\geq m$ and $\alpha_{2}\geq m$, then by
the identity (\ref{th3}) they can be represented as a linear
combination of the elements, in which $\alpha_{1}<m$. Thus any
monomial will be expressed through such monomials in which either
only $\alpha_{1}\geq m$ or $\alpha_{2}\geq m$. The number of such
monomials is bounded by the constant $2m$ independent of $n$.

In the general case the space $Q_{\lambda_{1}}$ will generate by
elements, in which only one $\alpha_{i}$ is not less then $tm$.
Note, that the general number of such elements is bounded by a
constant independent of $n$.

Let $i$ will be a smallest index for which $\alpha_{i}\geq tm$.
consider the corresponding element:
$$g_{s}=yY^{\alpha_{1}}x_{t_{1}}Y^{\alpha_{2}}\dots
Y^{\alpha_{i}}x_{t_{i}}Y^{\alpha{i+1}}\dots
Y^{\alpha_{s}}x_{t_{s}}Y^{\alpha_{s+1}}.$$ If $\alpha_{i+1}\geq tm$,
then the identity (\ref{th3}) allows to bring the element $g_{s}$ to
a linear combination of words, that are lexicographically
less.lexicographically smaller. If $\alpha_{i+1}<tm$, then modulo
words, lexicographically smaller, the element $g_{s}$ can be written
as $$Y^{\alpha_{1}}x_{t_{1}}\dots Y^{\alpha_{i}-\alpha_{i+1}-1}
\underbrace{(y(y\dots(yx_{t_{i}})\dots ))}_{\alpha_{i+1}+1}\cdot$$
$$\cdot
\underbrace{(y(y\dots(yx_{t_{i+1}})\dots))}_{\alpha_{i+1}}Y^{\alpha_{i+2}}\dots
Y^{\alpha_{s}}x_{t_{s}}Y^{\alpha_{s+1}}.$$ The Leibnitz identity
allows to bring the last element to the sum of terms, that are
lexicographically smaller, and term $$Y^{\alpha_{1}}x_{t_{1}}\dots
Y^{\alpha_{i}-\alpha_{i+1}-1} X'Y^{\alpha_{i+2}}\dots
Y^{\alpha_{s}}x_{t_{s}}Y^{\alpha_{s+1}},$$ where
$$x'=\underbrace{(y(y\dots(yx_{t_{i}})\dots ))}_{\alpha_{i+1}+1}
\underbrace{(y(y\dots(yx_{t_{i+1}})\dots))}_{\alpha_{i+1}}.$$ We
will contain element with fewer generators $x_{i_{r}}$ covered by
the induction assumption. The theorem is proved.

\section*{3. The basis of multilinear part of variety $\widetilde{\textbf{V}}_{3}$}
\section*{of Leibnitz algebras}

The variety $\widetilde{\textbf{V}}_{3}$ of Leibnitz algebras is an
equivalent to the well-known variety $\textbf{V}_{3}$ of Lie
algebras. Previously, in the article \cite{ALE-MSP} the growth of
this variety was designated, and in the article \cite{STV} --- its
multiplicities and colength.

Let $T=\Phi[t]$ be a ring of polynomial in the variable $t$.
Consider three-dimensional Heisenberg algebra $H$ with the besis
$\{a,b,c\}$ and multiplication $ba=-ab=c$, the product of the
remaining basis elements is zero. Well known and easy to verify that
the algebra $H$ is nilpotent of the class two Lie algebra. Transform
the polynomial ring $T$ in the right module of algebra $H$, in which
the basis elements of algebra $H$ act on the right on the polynomial
$f$ from $T$ follows:
$$fa=f', fb=tf, fc=f,$$ where $f'$ is a partial derivative
of a polynomial $f$ in the variable $t$. Consider the direct sum of
vector spaces $H$ and $T$ with multiplication by the rule:
$$(x+f)(y+g)=xy+fy,$$ where $x,y$ are from $H$; $f,g$ are from $T$. Denote it by the symbol
$\widetilde{H}$. Direct verification shows that $\widetilde{H}$ is
an algebra of Leibnitz. The algebra $\widetilde{H}$ is the Leibnitz
algebra, satisfies to the identity $x(y(zt))\equiv 0$ and generates
the variety $\widetilde{\textbf{V}}_{3}$ of Leibnitz algebras.

\textbf{Lemma.} \emph{The variety $\widetilde{\textbf{V}}_{3}$
satisfies to the identities:}
\begin {equation}
\label {3N4} x(y(zt))\equiv 0,
\end {equation}
\begin {equation}
\label {slova5}
x_{0}A\overline{x}_{1}B\overline{x}_{2}C\overline{x}_{3}D\overline{x}_{4}\equiv
0,
\end {equation}
\begin {equation}
\label {t-vo6} x_{0}(x_{1}x_{4})(x_{2}x_{3})\equiv
x_{0}(x_{1}x_{2})(x_{3}x_{4})+x_{0}(x_{1}x_{3})(x_{2}x_{4}),
\end {equation}
\emph{where $A,B,C,D$ are some words from generators.}

\textbf{Proof.} The truth of identities (\ref{3N4}) and
(\ref{slova5}) verified by arbitrary replacement generators by
elements of algebra $\widetilde{H}$ and was showed in the paper
\cite{ALE-MSP}. Consider the following special form of the second
identity:
$x_{0}\overline{x}_{1}\overline{x}_{2}\overline{x}_{3}\overline{x}_{4}\equiv
0.$ Presenting it as a sum and using the identity $xyz-xzy\equiv
x(yz)$, we obtain:
$2x_{0}(x_{1}x_{2})(x_{3}x_{4})-2x_{0}(x_{1}x_{3})(x_{2}x_{4})+2x_{0}(x_{1}x_{4})(x_{2}x_{3})\equiv
0.$ Dividing this identity by 2 and moving the second term to the
right, we obtain the identity (\ref{t-vo6}). The lemma is proved.

\textbf{Theorem 2.} \emph{The set elements of form}
$$\theta(i,i_{1},...,i_{m},j_{1},...,j_{m})=x_{i}(x_{i_{1}}x_{j_{1}})(x_{i_{2}}x_{j_{2}})...(x_{i_{m}}x_{j_{m}})x_{k_{1}}x_{k_{2}}...x_{k_{n-2m-1}},$$
\emph{where $i_{s}<j_{s}, s=1,2,...,m, i_{1}<i_{2}<...<i_{m},
j_{1}<j_{2}<...<j_{m}, k_{1}<k_{2}<...<k-2m-1$, is a basis of the
space $P_{n}(\widetilde{V}_{3})$.}

\textbf{Proof.} Consider an arbitrary element of the space
$P_{n}(\widetilde{\textbf{V}}_{3})$. Using corollary $xy(zt)\equiv
x(zt)y$ from the Leibnitz identity and identity (\ref{3N4}), move
the all pairs as far right as possible.

We order the elements obtained using the lexicographic ordering of
lines $(k_{1},k_{2},\dots,k_{n-2m-1})$. Let the considering element
has a form:
$$x_{i}(x_{i_{1}}x_{j_{1}})(x_{i_{2}}x_{j_{2}})...(x_{i_{m}}x_{j_{m}})x_{k_{1}}...x_{k_{s}}x_{k_{s+1}}...x_{k_{n-2m-1}}$$
and $k_{s}>k_{s+1}$. Using the identity $xyz\equiv xzy+x(yz)$ we can
write this element as a sum
$$x_{i}(x_{i_{1}}x_{j_{1}})(x_{i_{2}}x_{j_{2}})...(x_{i_{m}}x_{j_{m}})x_{k_{1}}...x_{k_{s+1}}x_{k_{s}}...x_{k-2m-1}+$$
$$+x_{i}(x_{i_{1}}x_{j_{1}})(x_{i_{2}}x_{j_{2}})...(x_{i_{m}}x_{j_{m}})x_{k_{1}}...x_{k_{s-1}}(x_{k_{s}}x_{k_{s+1}})x_{k_{s+2}}...x_{k-2m-1},$$
where the first term is lexicographically less, than parent element,
and the second term has fewer number of single elements. Applying
the same method to the resulting term, we eventually present our
original element as a sum of terms, in that
$k_{1}<k_{2}<...<k_{n-2m-1}$.

Consider an arbitrary element, in which indexes of single elements
are ordered. We choose the two lowest index in the considered
element and redenote them through  $1'$ and $2'$ relatively. We
introduce the lexicographic order on lines
$(j_{1},j_{2},\dots,j_{m})$. Using also the induction on the number
of brackets, we will prove, that all received elements can be
represented as a linear combination of elements
$\theta(i,i_{1},...,i_{m},j_{1},...,j_{m})$. The corollary
$x(yz)\equiv -x(zy)$ of Leibnitz identity allows to order the
indexes of elements in couples, and the identity $xy(zt)\equiv
x(zt)y$ allows to order the brackets by the indexes of first
elements. According to these identities, the element can be written
either in the form
$x_{i}(x_{1'}x_{2'})(x_{i_{1}}x_{j_{1}})...(x_{i_{m-1}}x_{j_{m-1}})x_{k_{1}}x_{k_{2}}...x_{k_{n-2m-1}},$
either in the form
$x_{i}(x_{1'}x_{j_{1}})(x_{2'}x_{j_{2}})...(x_{i_{m-2}}x_{j_{m}})x_{k_{1}}x_{k_{2}}...x_{k_{n-2m-1}}.$
In the first case we can consider the ordering on $m-1$ brackets
that runs by induction. In the second case we apply the identity
(\ref{t-vo6}) and obtain:
$x_{i}(x_{1'}x_{2'})(x_{j_{2}}x_{j_{1}})...(x_{i_{m-2}}x_{j_{m}})x_{k_{1}}x_{k_{2}}...x_{k_{n-2m-1}}+$
$x_{i}(x_{1'}x_{j_{2}})(x_{2'}x_{j_{1}})...(x_{i_{m-2}}x_{j_{m}})x_{k_{1}}x_{k_{2}}...x_{k_{n-2m-1}},$
where for the first term can be again apply the induction
hypothesis, and the second term is lexicographically less.
Therefore, any element of the space $P_{n}(\widetilde{V}_{3})$ can
be written as a linear combination of elements
$\theta(i,i_{1},...,i_{m},j_{1},...,j_{m})$ modulo
$Id(\widetilde{V}_{3})$.

We now prove, that the elements
$\theta(i,i_{1},...,i_{m},j_{1},...,j_{m})$ are linearly independent
modulo $Id(\widetilde{V}_{3})$. Consider the linear combination of
these elements:
$$\sum_{(i,i_{1},...,i_{m}.j_{1},...,j_{m})}\alpha(i,i_{1},...,i_{m},j_{1},...,j_{m})\theta(i,i_{1},...,i_{m},j_{1},...,j_{m})=0$$
and show that all coefficient
$\alpha(i,i_{1},...,i_{m},j_{1},...,j_{m})$ are zero. Assume the
contrary.

Choose \ an \ element \
$\theta(i^{*},i^{*}_{1},...,i^{*}_{m},j^{*}_{1},...,j^{*}_{m})$ \
with \ non-zero \ coefficient\
$\alpha(i^{*},i^{*}_{1},...,i^{*}_{m},j^{*}_{1},...,j^{*}_{m})$ such
that the number of commutators $m$ in it is the least and the index
$j^{*}_{1}$ of element in the second position in the first bracket
is the largest. Since each element is uniquely determined by the
number $m$, by the element $x_{i}$ and by the sample
$(i,i_{1},...,i_{m},j_{1},...,j_{m})$, then in the selected element
$\theta(i^{*},i^{*}_{1},...,i^{*}_{m},j^{*}_{1},...,j^{*}_{m})$
these rates are fixed. We replace its generators on the basis
elements of algebra $\widetilde{H}$ as follows: $x_{i^{*}}=f,
x_{i^{*}_{s}}=a, x_{j^{*}_{s}}=b, s=1,...,m$, the rest generators we
replace on element $c$. After this substitution all elements, which
are different from the chosen, will be equal to zero: if the element
$x_{i^{*}}$ will be replace on the basis element of Heisenberg
algebra, than this element will be zero (because $xf\equiv 0$ for
any $x$ from $\widetilde{H}$); if the element will have more then
$m$ commutators, then it will be also zero (since the element $c$
from the center of algebra fall into the commutator); a similar
situation arises, if the element will have $m$ commutators but the
sample will be different from the fixed. Indeed, there are two kinds
of elements containing $m$ commutators at a fixed sample
$(i^{*},i^{*}_{1},...,i^{*}_{m},j^{*}_{1},...,j^{*}_{m})$: these are
the elements that contain in the second position and the first
bracket $x_{j^{*}_{1}}$, and elements that contain in the second
position and the first bracket $x_{i^{*}_{s}}$, where
$i^{*}_{s}<j^{*}_{1}$ ($s=2,...,m$).

All elements of the second kind are zero, as by described
substitution the first bracket will be equal to $(aa)$. In one of
brackets of the elements of the first type fall generators
$x_{i^{*}_{s}}$ and $x_{i^{*}_{t}}$. As a result of described
substitution this bracket also resets the element. Thus we obtained,
that if $f\neq 0$, then
$\sum_{(i,i_{1},...,i_{m},j_{1},...,j_{m})}\alpha(i,i_{1},...,i_{m},j_{1},...,j_{m})\theta(i,i_{1},...,i_{m},j_{1},...,j_{m})=0.$
Consequently, contrary to the assumption
$\alpha(i^{*},i^{*}_{1},...,i^{*}_{m},j^{*}_{1},...,j^{*}_{m})$ is
zero. The theorem is proved.

Note that in the proof of the theorem we used only the Leibnitz
identity and corollaries from it, the identity $x(y(zt))\equiv 0$,
and lastly the identity $(3)$. Consequently, any identity that runs
in the variety $\widetilde{\textbf{V}}_{3}$, is a corollary from
these identities. Hence we obtain the following assertion.

\textbf{Corollary.} \emph{The identities (\ref{3N4}) and
(\ref{t-vo6}) form a basis of identities of variety
$\widetilde{\textbf{V}}_{3}$.}

The authors thank S.P. Mishchenko for useful advice and attention to
this work.

\end{document}